\newenvironment{E}{\begin{equation}}{\end{equation}}
\def\proof{\noindent{\bf Proof: }}
\def\qed{ \hskip 20pt{\vrule height7pt width6pt depth0pt}\hfil}
\def\forb{{\hbox{forb}}}

\def\0{{\bf 0}}
\def\1{{\bf 1}}

\def\s{{\hbox{supp}}}
\def\Av{{\mathrm{Avoid}}}
\newcommand{\linelessfrac}[2]{\genfrac{}{}{0pt}{}{#1}{#2}}
\newcommand{\ncols}[1]{\| #1 \|}
\newcommand{\rf}[1]{(\ref{#1})}
\newcommand{\trf}[1]{Theorem~\ref{#1}}
\newcommand{\lrf}[1]{Lemma~\ref{#1}}

\newcommand{\corf}[1]{Conjecture~\ref{#1}}

\newcommand{\srf}[1]{Section~\ref{#1}}

\documentclass[12pt]{article}
\newtheorem{thm}{Theorem}[section]
\newtheorem{lemma}[thm]{Lemma}

\newtheorem{conj}[thm]{Conjecture}
\newtheorem{defn}[thm]{Definition}

\renewcommand{\l}{\ell}
\usepackage{amssymb}
\usepackage{amsmath}
\title{ Repeated columns and an old chestnut  }
\author{R.P. Anstee\thanks{Research supported in part by
NSERC, work done while visiting the second author at USC.}
\\Mathematics Department\\The University of British Columbia\\Vancouver,
B.C. Canada V6T 1Z2\\ \\
\and Linyuan Lu\thanks{
This author was supported in part by NSF
grant  DMS 1000475. }
 \\Mathematics Department\\The University of South Carolina\\Columbia, SC, USA  \\
\\\mbox{\ }}

\begin{document}
\maketitle
\begin{abstract}
Let $t\ge 1$ be a given integer. Let ${\cal F}$ be a family of  subsets of $[m]=\{1,2,\ldots ,m\}$.  Assume that for every pair of disjoint sets $S,T\subset [m]$  with $|S|=|T|=k$, 
there  do not exist   $2t$ sets in ${\cal F}$ where $t$ subsets of ${\cal F}$ contain $S$ and are disjoint from $T$ and $t$ subsets of ${\cal F}$ 
contain $T$ and are  disjoint from $S$. We show that $|{\cal F}|$ is $O(m^{k})$. 

Our main new ingredient is allowing, during the inductive proof, multisets of subsets of $[m]$ where the multiplicity of a given set is bounded by $t-1$.  We use a strong stability result of Anstee and Keevash. This is further evidence for  a conjecture of Anstee and Sali.  These problems can be stated in the language of matrices  Let $t\cdot M$ denote $t$ copies of the matrix $M$ concatenated together.  We have established the conjecture for those configurations $t\cdot F$ for any $k\times 2$ (0,1)-matrix $F$.

\end{abstract}

\section{Introduction}

We will be considering a problem in extremal hypergraphs that can be phrased as how many edges  a hypergraph  on $m$ vertices can have when there is a forbidden subhypergraph. There are a variety of ways to define this problem (we could, but do not, restrict to (simple) $k$-uniform hypergraphs).  We can encode a hypergraph on $m$ vertices
as an $m$-rowed (0,1)-matrix where the $i$th column is the incidence vector of  the $i$th hyperedge. A hypergraph is \emph{simple} if there are no repeated edges. We define 
 a matrix to be \emph{simple} if it is a (0,1)-matrix with no repeated columns.
We will use the language of matrices in this paper.

Let $M$ be an $m$-rowed (0,1)-matrix. Some notation about repeated columns is needed.
For an $m\times 1$ (0,1)-column $\alpha$, we define $\mu(\alpha,M)$ as
the multiplicity of column $\alpha$ in a matrix $M$. We  consider
matrices of bounded column multiplicity.  We define a matrix $A$ to be
$t$-\emph{simple} if it is a (0,1)-matrix and every column $\alpha$ of
$A$ has $\mu(\alpha,A)\le t$. Simple matrices  are 1-\emph{simple}.
For a given matrix $M$, let $\s(M)$ denote the maximal simple
$m$-rowed submatrix of $M$, so that if $\mu(\alpha,M)\ge 1$ then
$\mu(\alpha,\s(M))=1$. The matrices below are a $3$-simple matrix $M$
and its support $\s(M)$. 
$$M=
   \left[\begin{array}[c]{cccccc}
      0 & 1 & 0 & 1 &1 & 0\\ 
      1 & 0 & 1 &1 & 0 & 1\\
    \end{array} \right],\qquad
    \s(M)=\left[   \begin{array}[c]{ccc}
      0 & 1 & 1\\
      1 & 0 & 1 \\
    \end{array}\\
    \right]$$

For two  $(0,1)$-matrices $F$ and $A$, we say that $F$ is a \emph{configuration} in $A$, and write $F \prec A$ if there is a row and column permutation of $F$ which is a submatrix of $A$.  Let ${\cal F}$ denote a finite set of (0,1)-matrices. Let $\Av(m,{\cal F},t)$ denote all $m$-rowed $t$-simple matrices $A$ for which $F\not\prec A$ for all $F\in{\cal F}$. We are most interested in cases with $|{\cal F}|=1$ \cite{survey}. We do not require any $F\in{\cal F}$ to be simple which is quite different from usual forbidden subhypergraph problems.
Our extremal function of interest is
$$\forb(m,{\cal F})=\max_A\{\ncols{A}\,:\,A\in\Av(m,{\cal F},1)\}.$$
We find it helpful to also define
$$\forb(m,{\cal F},t)=\max_A\{\ncols{A}\,:\,A\in\Av(m,{\cal F},t)\}.$$
If $A\in\Av(m,{\cal F},t)$ then $\s(A)\in\Av(m,{\cal F},1)$ and 
$\ncols{A}\le t\cdot \ncols{\s(A)}$.  We obtain 
\begin{E}\forb(m,{\cal F})\le \forb(m,{\cal F},t)\le t\cdot\forb(m,{\cal F}),\label{asymptotics}\end{E}
 so that the asymptotic growth  of $\forb(m,{\cal F})$ is the same as that of $\forb(m,{\cal F},t)$.

We have an important conjecture about $\forb(m,F)$.  We use the notation
$[M\,|\,N]$  to denote the matrix obtained from concatenating the two
matrices $M$ and $N$. We use the notation $k\cdot M$ to denote the
matrix $[M|M|\cdots |M]$ consisting of $k$ copies of $M$ concatenated
together.
Let $I_k$ denote the $k\times k$ \emph{identity} matrix and let $I_k^c$ denote the (0,1)-complement of $I_k$. Let $T_k$ denote the  $k\times k$ 
\emph{triangular} (0,1)-matrix with the $(i,j)$ entry being 1 if and only if $i\le j$. For an $m_1\times n_1$ matrix $X$ and an $m_2\times n_2$ matrix $Y$, we define the $2$-fold product $X\times Y$ as the $(m_1+m_2)\times n_1n_2$ matrix   each column consisting of a column of $X$ placed on a column of $Y$ and this is done in all possible ways. This  extends to $p$-fold products.

\begin{defn} Let $X(F)$ be the smallest $p$
so that $F\prec A_1\times A_2\times \cdots \times A_p$ for every choice of $A_i$ as either
$I_{m/p}$, $I_{m/p}^{c}$ or $T_{m/p}$. \end{defn}
Alternatively, assuming $F\not\prec I$ or $F\not\prec I^c$ or $F\not\prec T$, then $X(F)-1$ is the largest choice of $p$ so that $F\not\prec A_1\times A_2\times \cdots \times A_p$ for some choices of $A_i$ as either
$I_{m/p}$, $I_{m/p}^{c}$ or $T_{m/p}$. We note that if  $A_1\times A_2\times \cdots \times A_p\in\Av(m,F)$, then $\forb(m,F)$ is $\Omega(m^p)$.

Details are in \cite{survey}. We are assuming $m$ is
large and divisible by $p$, in particular that $m\ge (k+1)(k\l+1)$ so that $m/p\ge k\l+1$. Divisibility by $p$ does not affect the asymptotic growth, thus
$\forb(m,F)$ is $\Omega(m^{X(F)-1})$ using an appropriate $(X(F)-1)$-fold product. 

\begin{conj}\label{grand}\cite{AS05} Let $F$ be given. Then  
$\forb(m,F)= \Theta(m^{X(F)-1}).\qed$\end{conj}

The conjecture was known to be true for all $3$-rowed $F$ \cite{AS05}  and all $k\times 2$ $F$ \cite{AK}.  \srf{evidence} shows how \trf{chestnut} establishes the conjecture for  matrices $t\cdot F$ when $F$ is a $k\times 2$ matrix. 
It is of interest to generalize \corf{grand} to $\forb(m,{\cal F})$ where $|{\cal F}|>1$ but we know example of ${\cal F}$ where the conjecture fails.

We define $F_{e,f,g,h}$ as the $(e+f+g+h)\times 2$ matrix consisting of $e$ rows $[1\,1]$, $f$ rows $[1\,0]$, $g$ rows $[0\,1]$ and $h$ rows $[0\,0]$.
 Let $\1_e\0_f$ denote the $(e+f)\times 1$ vector of $e$ 1's on top of $f$ 0's so that $F_{e,f,g,h}=[\1_{e+f}\0_{g+h}\,|\,\1_{e}\0_{f}\1_{g}\0_{h}]$.
We let $\1_e$ denote the $e\times 1$ vector of $e$ 1's and $\0_f$
denote the $f\times 1$ vector of $f$ 0's.
 Our main result is the following which had foiled many previous attempts.
\begin{thm}\label{chestnut}Let $t\ge 2$ be given. Then $\forb(m,t\cdot
  F_{0,k,k,0})$ is $\Theta(m^k)$.\end{thm}

The forbidden configuration $t\cdot F_{0,k,k,0}$ in the language of sets,  consists
of two disjoint $k$-sets $S$, $T$, and a family of $t$ sets containing
$S$ but disjoint from $T$, and the other family of another $t$ sets
containing $T$ but disjoint from $S$. This theorem echoes our statement
in the abstract.

The result for $t=2$ and $k=2$ was proven in \cite{A90} and many details worked out for $t=2$ and $k>2$ by the first author and Peter Keevash. The extension for $t>2$, $k=2$ had been open since then \cite{survey}. The proof for $t>2$, $k=2$ is in \srf{induction}. The proof  for $t>2$, $k>2$ is in \srf{evidence}. Matrices $F_6(t),F_7(t)$ were given in \cite{survey}   as 4-rowed forbidden configurations (with some columns of multiplicity $t$) for which \corf{grand} predicts
$\forb(m,F_6(t))$ and $\forb(m,F_7(t))$ are $O(m^2)$.   Note that  $t\cdot F_{0,2,2,0}\prec F_6(t)$ and $t\cdot F_{0,2,2,0}\prec F_7(t)$ and so  \trf{chestnut} is a step towards these bounds which would establish \corf{grand} for all $4$-rowed $F$. 
Our proof  use a new induction given in \srf{induction} that considers $t$-simple matrices as well as a strong stability result \lrf{strongstability}. We offer some additional applications in \srf{additional}.

\section{ New Induction}\label{induction}

We  consider a new form of the standard induction for forbidden configurations \cite{survey}. Let $F$ be a matrix with maximum column multiplicity $t$.  Thus $F\prec t\cdot \s(F)$.
Let $A\in\Av(m,F,t-1)$. Assume $\ncols{A}=\forb(m,{\cal F},t-1)$. Given a row $r$ we permute rows and columns of $A$ to obtain

\begin{E}A=\begin{matrix}\hbox{ row }r\rightarrow\\ \\ \end{matrix}
\left[\begin{matrix}0\,0\cdots\, 0&1\,1\cdots \,1\\ G&H\end{matrix}\right].\label{rdecompbasic}\end{E}
Now $\mu(\alpha,G)\le t-1$ and $\mu(\alpha,H)\le t-1$.  For those $\alpha$ for which $\mu(\alpha,[G\,H])\ge t$, let $C$ be formed with
$\mu(\alpha,C)=\min\{\mu(\alpha,G),\mu(\alpha,H)\}$.  
We rewrite our decomposition of $A$ as follows:
\begin{E}A=\begin{matrix}\hbox{ row }r\rightarrow\\ \\ \end{matrix}
\left[\begin{matrix}0\,0\cdots\, 0&1\,1\cdots \,1\\ B\quad C&C\quad D\end{matrix}\right].\label{rdecompt}\end{E}
Then we deduce that $[BCD]$ and $C$ are both  $(t-1)$-simple. The former follows from
$\mu(\alpha,[B\,C\,D])=\mu(\alpha,G)+\mu(\alpha,H)- \min\{\mu(\alpha,G),\mu(\alpha,H)\}\le t-1$. We have that $F\not\prec [B\,C\,D]$ for $F\in{\cal F}$. 
Also for any $F'\prec C$ then  $[0\,1]\times F'\prec A$ so we define 
\begin{E}{\cal G}=\{F'\,:\,\hbox{for }F\in{\cal F},\,\, F\prec [0\,1]\times F'\hbox{ and }F\not\prec [0\,1]\times F''\hbox{ for all }F''\prec F',\,\,F''\ne F'\}.\label{inductivechildren}\end{E} 
Also
since each column $\alpha$ of $C$ has $\mu(\alpha,[G\,H])\ge t$, we deduce that $\s(F)\not\prec C$ for each $F\in{\cal F}$. Our induction on $m$ becomes:
$$\forb(m,{\cal F},t-1)=\ncols{A}=\ncols{[BCD]}+\ncols{C}\hskip 2.5in$$
\begin{E}\hskip 1in\le \forb(m-1,{\cal F},t-1)+ (t-1)\cdot\forb(m-1,{\cal G}\cup\{\s(F)\,:\,F\in{\cal F}\}).\label{newinductiont}\end{E}

\noindent{\bf Proof of \trf{chestnut} for $k=2$: } We will use induction on $m$ to show $\forb(m,t\cdot F_{0,2,2,0},t)$ is $O(m^2)$.  The maximum multiplicity of a column in 
$t\cdot F_{0,2,2,0}$ is $t$ and $F_{0,2,2,0}=\s(t\cdot F_{0,2,2,0})$.  Also $t\cdot F_{0,2,2,0}\prec [0\,1]\times (t\cdot F_{0,2,1,0})$. Let $A\in\Av(m,t\cdot F_{0,2,2,0},t-1)$ with $\ncols{A}=\forb(m,t\cdot F_{0,2,2,0},t-1)$. Apply \rf{newinductiont}.
We have 
$$\forb(m,t\cdot F_{0,2,2,0},t-1)=\ncols{A}=\ncols{[BCD]}+\ncols{C}\hskip 2.5in$$
$$\hskip 1in\le \forb(m-1,t\cdot F_{0,2,2,0},t-1)+(t-1)\cdot \forb(m-1,\{F_{0,2,2,0},t\cdot F_{0,2,1,0}\}).$$
  We apply
\lrf{linearbd} with induction on $m$ to deduce that $\forb(m,t\cdot F_{0,2,2,0},t-1)$ is $O(m^2)$. Then by \rf{asymptotics}, $\forb(m,t\cdot F_{0,2,2,0})$ is also $O(m^2)$.\qed

\vskip 10pt
\trf{chestnut} was proven for $t=k=2$ in \cite{A90} using induction in the spirit \rf{newinductiont} ($(t-1)$-simple matrices are simple) 
and 
\lrf{linearbd} for $t=2$.

\begin{lemma}We have that $\forb(m,\{F_{0,2,2,0},t\cdot F_{0,2,1,0}\})$ is $O(m)$. \label{linearbd}\end{lemma}
\proof
Let $A\in\Av(m,\{F_{0,2,2,0},t\cdot F_{0,2,1,0}\})$. Avoiding $F_{0,2,2,0}$ creates structure:
Let $X_i$ denote the columns of $A$ of column sum $i$.    Let $J_{a\times b}$ denote the 
$a\times b$ matrix of 1's and let $0_{a\times b}$ denote the $a\times 
b$ matrix of 0's.  Now $F_{0,2,2,0}\not\prec X_i$ and so for $\ncols{X_i}\ge 3$, we may deduce that there is a partition of the rows $[m]$ into 
$A_i\cup B_i\cup C_i$. Let $x_i=|X_i|$.  After suitable row and column permutations, we have $X_i$ as follows:
$$\hbox{type 1: }
X_i=\begin{array}{l@{}} A_i\{\\ B_i\{\\ C_i\{\\ \end{array}
\hskip-3pt
\left[\begin{array}{l}
I_{x_i}\\ J_{(i-1)\times x_i}\\ 0_{(m-x_i-i+1)\times x_i}\\
\end{array}\right]
\hbox{ or }
\hbox{type 2: }
X_i=\begin{array}{l@{}} A_i\{\\ B_i\{\\ C_i\{\\ \end{array}
\left[\begin{array}{l}
I_{x_i}^c\\ J_{(i-x_i+1)\times x_i}\\ 0_{(m-i-1)\times  x_i}\\
\end{array}\right]
.$$
We will say $i$ \emph{ is of type } $j$  ($j=1$ or $j=2$) if the columns 
of sum $i$ are of type $j$. These are the \emph{sunflowers} (for type 1) and \emph{inverse sunflowers} (type 2) of \cite{FFP} where for type 1 the petals are $A_i$ with center $B_i$.

Let $T(1)=\{i\,:\,i\hbox{ is of type 1 and }\ncols{X_i}\ge t+2\}$.
We wish to show for  that $B_{i}\subset B_{j}$ for $i,j\in T(1)$ and $i<j$. Assume $p\in B_{i}\backslash  B_{j}$. 
Given that $|B_{i}|<|B_{j}|$, there are two rows $r,s\in  B_{j}\backslash  B_{i}$. Then we find a copy of $t\cdot F_{0,2,1,0}$ in rows $p,r,s$
of $[X_{i}\,X_{j}]$ (we would not choose the possible column of $X_i$ that has a 1 in row $r$ and the   column of $X_i$ that has a 1 in row $s$), a contradiction showing no such $p$ exists and hence $B_{i}\subset B_{j}$.

We form a matrix $Y_1$ from those $X_i$ with $i\in T(1)$. We have $\ncols{Y_1}=\sum_{i\in T(1)} \ncols{X_i} =\sum_{i\in T(1)}|{A_{i}}|$. Assume $\sum_{i\in T(1)}|{A_{i}}|> (t+1)m$. 
Then there is some row $p$ and $(t+2)$-set  $\{s(1),s(2),\ldots, s(t+2)\}$ with $p\in A_i$ for all $i\in\{s(1),s(2),\ldots, s(t+2)\}$. Assume $s(1)<s(2)<\cdots<s(t+2)$. 
We have $B_{s(1)}\subset B_{s({2})}\subset\cdots\subset  B_{s({t+2})}$. 
We may choose $r,s\in B_{s({t+2})}\backslash B_{s({t})}$ so that $r,s\in A_{s(i)}\cup C_{s(i)}$ for $i=1,2,\ldots ,t$. We find  a copy of $t\cdot F_{0,2,1,0}$ in rows  $p,r,s$ as follows. We take one column from each  $X_{s(j)}$ for $j=1,2,\ldots ,t$ and $t$ columns from the $X_{s(t+2)}$.   We conclude that $\ncols{Y_1}\le (t+1)m$.  Similarly the matrix $Y_2$ formed from those $X_i$  such that $i$ is of type 2 and 
$\ncols{X_i}\ge t+2$ has $\ncols{Y_2}\le (t+1)m$.  Now $Y_1$ and $Y_2$ represent all columns of $A$ with the exception of columns of sum $i$ with $\ncols{X_i}\le t+1$ and so we conclude $\ncols{A}\le \ncols{Y_1}+\ncols{Y_2}+(t+1)(m-1)+2$. Thus $\ncols{A}$ is $O(m)$.
\qed

\section{More evidence for the Conjecture}\label{evidence}
 This section first explores the \corf{grand} for $t\cdot F$ when $F$ is $k\times 2$.  The section concludes with the proof of \trf{chestnut} for $k>2$.
The following verifies \corf{grand} for all $k\times 2$ $F$.
Note that any $k\times 2$ matrix $F$ can be written as $F_{a,b,c,d}$
($b\geq c$) under proper row and column permutations. Since
$\forb(m,F)$ is invariant under taking $(0,1)$-complement, we can
further assume $a\geq d$. The case of $t=1$ was solved in \cite{AK} by
the following theorem.

\begin{thm}\label{twocol}\cite{AK}
Suppose $a \geq d$ and $b \geq c$. Then
$\forb(m,F_{a,b,c,d})$ is $\Theta(m^{a+b-1})$ if either $b>c$ or $a,b\ge 1$.
Also $\forb(m,F_{a,0,0,d})$ is $\Theta(m^a)$ and
$\forb(m,F_{0,b,b,0})$ is $\Theta(m^{b})$. \qed
\end{thm}

Note that \corf{grand} is verified if there is a product construction avoiding $F$ yielding the same asymptotic growth as an upper bound on $\forb(m,F)$. 
The $k$-fold product $I_{m/k}\times I_{m/k}\times\cdots\times I_{m/k}\in\Av(m,t\cdot F_{0,k,k,0})$ has $\Theta(m^k)$ columns. Thus \trf{chestnut} verifies the conjecture for $t\cdot F_{0,k,k,0}$.  The following results  verify the conjecture for $t\cdot F$ for the remaining $k\times 2$ $F$. 
 
\begin{thm} For $b>c$ or $a,b\ge 1$ then $\forb(m,t\cdot F_{a,b,c,d})$ is $\Theta(m^{a+b})$.\end{thm}
\proof  The upper bound follows from $\forb(m, F_{a,b,c,d})$ being $\Theta(m^{a+b-1})$ and then applying \lrf{ttimes}.  
The lower bound follows from $2\cdot\1_{a+b}\prec t\cdot F_{a,b,c,d}$  so that the $(a+b)$-fold product 
$I_{m/(a+b)}\times I_{m/(a+b)}\times\cdots\times I_{m/(a+b)}\in\Av(m,F_{a,b,c,d})$ and hence $\forb(m,t\cdot F_{a,b,c,d})$ is $\Omega(m^{a+b})$.\qed

\begin{thm} Let $a\ge d$ be given. Then $\forb(m,t\cdot  F_{a,0,0,d})$ is $\Theta(m^a)$.\end{thm}
\proof This follows using \lrf{10} repeatedly and also $\forb(m,t\cdot  F_{a,0,0,0})$ is $O(m^a)$ using \trf{general}. 
The $a$-fold product $I_{m/a}\times I_{m/a}\times\cdots\times I_{m/a}\in \Av(m,t\cdot  F_{a,0,0,d})$. \qed

\vskip 10pt
The following result can be found in the survey on forbidden configurations \cite{survey}

\begin{lemma} Assume $\forb(m,F)$ is $O(m^{\ell})$. Then $\forb(m,\left[\linelessfrac{1}{0}\right]\times F)$ is $O(m^{\ell+1})$.\label{10}\end{lemma}

Here is the summary of results on $\forb(m, t\cdot F_{a,b,c,d})$
($a\geq d$ and $b\geq c$),
which verify \corf{grand} for all $k\times 2$ $F$.

\begin{table}[htb]
  \centering
  \begin{tabular}[c]{|c|c|c|c|c|}
 \hline  \hline
$t$ & Configuration & result & reference & Lower bound construction \\
\hline
 & \parbox{2.8cm}{$F_{a,b,c,d}$ ($b>c$ \mbox{ or } $a,b\geq 1$)}  
& $\Theta(m^{a+b-1})$ & \cite{AK} &
$\overbrace
{I\times I \times \cdots I \times I}^{a+b-1}$\\
\cline{2-5}
$t=1$ & \parbox{2.8cm}{$F_{a,0,0,d}$} & $\Theta(m^{a})$ & \cite{AK} &
$\overbrace
{I\times I \times \cdots I \times I}^{a}$\\
\cline{2-5}
& \parbox{2.8cm}{$F_{0,b,b,0}$} & $\Theta(m^{b})$ & \cite{AK} &
$\overbrace
{I\times I \times \cdots I \times I}^{b-1}\times T$\\
\hline
\hline
 & \parbox{2.8cm}{$t\cdot F_{a,b,c,d}$ ($b>c$ \mbox{ or } $a,b\geq 1$)}  
& $\Theta(m^{a+b})$ & \lrf{ttimes} &
$\overbrace
{I\times I \times \cdots I \times I}^{a+b}$\\
\cline{2-5}
$t\geq 2$ & \parbox{2.8cm}{$t\cdot F_{a,0,0,d}$} & $\Theta(m^{a})$ & \lrf{ttimes}  &
$\overbrace
{I\times I \times \cdots I \times I}^{a}$\\
\cline{2-5}
& \parbox{2.8cm}{$t\cdot F_{0,b,b,0}$} & $\Theta(m^{b})$ & Theorem
\ref{chestnut}  &
$\overbrace
{I\times I \times \cdots I \times I}^{b-1}\times T$\\
\hline

  \end{tabular}
  \caption{All cases of $\forb(m, t\cdot F_{a,b,c,d})$ with $a\geq d$
    and $b\geq c$.}
  \label{tab:1}
\end{table}

\vskip 5pt  
We note that the bound for $\forb(m,t\cdot F_{a,0,0,d})$ can be readily established by a pigeonhole argument.
We return to \trf{chestnut} and first obtain some useful lemmas.
 Let $X_i\in\Av(m,F_{0,k,k,0})$ with all column sums $i$. We define $X_i$ to be of \emph{type} $(a,b)$ if $a,b\ge 0$ are integers with $a+b=k-1$ and there is a partition
$C_i\cup D_i=[m]$ with $|D_i|+a-b=i$ such that  any column $\alpha$ of $X_i$ has  exactly $a$ 1's in rows $C_i$ and  exactly $b$ 0's in rows $D_i$.  We are able to use this structure in view of the following `strong stability' result:

\begin{lemma}\cite{AK} Let $Y_i\in\Av(m,F_{0,k,k,0})$ with all column sums $i$. Assume $\ncols{Y_i}\ge (6(k-1))^{5k+2}m^{k-2}$. Then 
there is an $m$-rowed submatrix $X_i$ of $Y_i$ and a pair of integers $a,b\ge 0$ with $a+b=k-1$ such that $X_i$ is of type $(a,b)$ and   where $\ncols{Y_i}-\ncols{X_i}\le m^{k-3}$.
\label{strongstability}\end{lemma}

\begin{lemma}\label{CD}Let $X_i\in\Av(m,F_{0,k,k,0})$ have all columns of sum $i$ and assume $X_i$ is of type $(a,b)$ with $a,b\ge 1$ with $a+b=k-1$. Let  $C_i\cup D_i=[m]$ be the associated partition of the rows.  We form a bipartite graph $G_i=(V_i,E_i)$ with $V_i=\binom{C_i}{a}\cup\binom{D_i}{b}$ where we have  $(C,D)\in E_i$ if there is a column of $X_i$ with 
$a$ 1's in rows $C$ and $D_i\backslash D$ and  $b$ 0's in rows $D$ and $C_i\backslash C$. 
Assume $|E_i|\ge 2km^{k-2}$.Then there is subgraph $G_i'=(V_i',E_i')$ of $G_i$ with
$|E'|\ge \frac{1}{2} |E_i|$ such that for every  pair $C\in\binom{C_i}{a}$ and $D\in\binom{D_i}{b}$ with
$(C,D)\in E'$ we have
\begin{E}d_{G_i'}(C)\ge (b+1/2)m^{b-1},\qquad d_{G_i'}(D)\ge (a+1/2)m^{a-1}.\hfill \label{mindCmindD}\end{E}
 \end{lemma}
  \proof    Simply delete vertices $C\in\binom{C_i}{a}$ with $d_G(C)< (b+1/2)m^{b-1}$ and vertices $D\in\binom{D_i}{b}$ with $d_G(D)< (a+1/2)m^{a-1}$ and continue deleting vertices until conditions
  \rf{mindCmindD} are satisfied for any remaining vertices of $G'$.  This will delete a maximum of
  $(b+1/2)m^{b-1}\binom{|C_i|}{a}+(a+1/2)m^{a-1}\binom{|D_i|}{b}<km^{k-2}$ edges which deletes less than half the edges of $G$.
  \qed

\begin{lemma} Let $k$ be given. Then $\forb(m,\{ F_{0,k,k,0},t\cdot F_{0,k,k-1,0}\})$ is $O(m^{k-1})$.  \label{smallbd}\end{lemma}
\proof  
Let $A\in\Av(m,\{ F_{0,k,k,0},t\cdot F_{0,k,k-1,0}\})$. Let $Y_i$ denote the columns of $A$ of column sum $i$.
For all $i$ for which $|Y_i|< (6(k-1))^{5k+2}m^{k-2}$, delete the columns of $Y_i$ from $A$. This may delete $ (6(k-1))^{5k+2}m^{k-1}$ columns.
For $i$ with $|Y_i|\ge (6(k-1))^{5k+2}m^{k-2}$,   apply \lrf{strongstability} and obtain $X_i$ with $|X_i|\ge  (6(k-1))^{5k+2}m^{k-2}-m^{k-3}$. 

  We consider a choice $a,b$ with $a+b=k-1$. Let $T(a,b)=\{i\,:\,X_i \hbox{ is of type }(a,b)\}$. 
We will show that $\sum_{i\in T(a,b)}|X_i|\le (tk)m^{k-1}$.

\noindent {\bf Case 1}.  $a,b\ge 1$.  
  
  Create $G_i$ as described in \lrf{CD} to obtain $G_i'$ for each $i\in T(a,b)$.  Now if $\sum_{i\in T(a,b)}|E_i'|> (t+1)m^{a+b}$, then there will be some  edge $(C,D) \in E_i'$ for at least $t+2$ choices $i\in T(a,b)$. Let those choices be $s(1),s(2),\ldots ,s(t+2)$ where $s(1)<s(2)<\cdots<s(t+2)$.   We wish to show that $X_{s(i)}$ has $t\cdot F_{0,k-1,0,0}$  on rows $C\cup D$.
  $$
  \begin{array}{c}
  \begin{array}{c@{}ccc}
 \hbox{rows }C\Biggl\{
  &
  \begin{array}{c}
  1\\
  1\\
  1\\
  \end{array}
 &
  \overbrace{\begin{array}{c}
  1\,1\cdots 1\\
  1\,1\cdots 1\\
  1\,1\cdots 1\\
  \end{array}}^{t}
   &
   \overbrace{\begin{array}{c} 
  0\,0\cdots 0\\
  0\,0\cdots 0\\
  0\,0\cdots 0\\
  \end{array}}^{t}
  \\
  \end{array}\\
   \begin{array}{c@{}ccc}
 \hbox{rows }D\biggl\{
  &
  \begin{array}{c}
  0\\
  0\\
  \end{array}
 &
  \begin{array}{c}
  1\,1\cdots 1\\
  1\,1\cdots 1\\
  \end{array}
   &
  \begin{array}{c} 
  0\,0\cdots 0\\
  0\,0\cdots 0\\
  \end{array}
  \\
  \end{array}\\
  \end{array}
  $$
  For a given set $D\in\binom{D_{s(i)}}{b}$, we compute
  $|\{H\in \binom{D_{s(i)}}{b}\,:\,H\cap D\ne\emptyset\}|\le\sum_{j=1}^{b}\binom{b}{j}\binom{D_{s(i)}\backslash D}{b-j}<bm^{b-1}$. 
  
  Now if 
  $d_{G'}(C)\ge (b+1/2)m^{b-1}$ and $(C,D)\in E_{s(i)}'$ then there are at least $t$ edges $(C,H)\in E_{s(i)}'$ with $H\cap D=\emptyset$. We are using
$(b+1/2)m^{b-1}>bm^{b-1} +t+2$ which is true for $m$ large enough and so asymptotics are unaffected.
  Thus we have 
  $t$ columns of  $X_{s(1)}$ with $\1_{k-1}$ on rows $C\cup D$ and, because these columns have  a 1's on rows $C\subseteq C_{s(1)}$, these columns are 0's on the remaining rows of $C_{s(1)}\backslash C$.  
  
  Similarly, because $d_{G_i'}(D)\ge (a+1/2)m^{a-1}$ there will be $t+2$ edges $(K,D)\in E_{s(i)}$ with $K\cap C=\emptyset$ and so there are  $t$ columns of $X_{s(t+2)}$ with  $\0_{k-1}$ on rows $C\cup D$ and, because these columns have 0's on rows $D$, these columns are 1's on rows of $D_{s(t+2)}\backslash D$.   
  
  We choose $k$ rows in
  $Z=D_{s(t+2)}\backslash D_{s(1)}$ so that $Z\subseteq C_{s(1)}$.   We deduce that in the chosen $t$ columns of  $X_{s(1)}$ we have $\0_k$ in rows $Z$ since 
  $Z\subseteq C_{s(1)}\backslash C$ and the columns have $\1_{k-1}$ in rows $C\cup D$. In the chosen $t$ columns  of $X_{s(t+2)}$ we have $\1_k$ in rows $Z$ since $Z\subset D_{s(t+2)}\backslash D$ and the columns have $\0_{k-1}$ in rows $C\cup D$.  This yields $t\cdot F_{0,k,k-1,0}$, a contradiction. Thus $\sum_{i\in Type(a,b)}|E_i'|\le (t+1)m^{k-1}$. 
This concludes Case 1.

\vskip 5pt

  \noindent {\bf Case 2}.  $a=k-1, b=0$ or $a=0, b=k-1$.   
 
    We proceed similarly.  We need only consider $a=k-1, b=0$ since the case $a=0, b=k-1$ is just the (0,1)-complement. For $i\in T(k-1,0)$, $X_i$ has partition $C_i\cup D_i=[m]$ and columns of $X_i$ have 1's on exactly $k-1$ rows of $C_i$ and all 1's on rows $D_i$. Assume   $\sum_{i\in T(k-1,0)}|X_i|\ge (tk)m^{k-1}$.  Then there are  $tk$ choices $s(1),s(2),\ldots ,s(tk)\in T(k-1,0)$ where $s(1)<s(2)<\cdots<s(tk)$ such that, for some $C\in \binom{C_s(i)}{k-1}$, each $X_{s(i)}$ has  a column with 1's in rows $C\cup D_{s(i)}$ and 0's in rows $C_{s(i)}\backslash C$.  We wish to find $t\cdot F_{0,k-1,0,0}$ in $A$ in rows $C$ as follows using one column from each of $X_{s(i)}$ for $i=1,2,\ldots ,t$ and $t$ columns from
    $X_{s(tk)}$.
$$
\begin{array}{ccccc@{}c}
\begin{array}{c}
\hbox{rows }C\,\left\{\begin{array}{@{}c@{}} \\ \\ \\ \end{array}\right.\\
\\
\end{array}
\begin{array}{c}
1\\
1\\
1\\
 X_{s(1)}\\
\end{array}
\begin{array}{c}
1\\
1\\
1\\
 X_{s(2)}\\
\end{array}
\begin{array}{c}
\\
\ldots\\
\\
\\
\end{array}
\begin{array}{c}
1\\
1\\
1\\
 X_{s(t)}\\
\end{array}
\begin{array}{c}
1\\
1\\
1\\
X_{s(tk)}\\
\end{array}&
\overbrace{\begin{array}{c}
0\,0\cdots 0\\
0\,0\cdots 0\\
0\,0\cdots 0\\
 X_{s(tk)}\\
\end{array}}^{t}\\
\end{array}
$$
     Given our choice $C\in \binom{C_{s(tk)}}{k-1}$, we compute that  $|\{K\in \binom{C_{s(kt)}}{k-1}\,:\, K\cap C\ne\emptyset\}|< km^{k-2}$. Thus with $|X_{s(kt)}|\ge km^{k-2}$,  there will be $t$ choices $K_1,K_2,\ldots ,K_t$ disjoint from $C$ and hence  one 
    column of $X_{s(kt)}$ for each $i=1,2,\ldots ,t$ with $\1_{k-1}$ on rows of $K_i\subseteq C_{s(kt)}\backslash C$
    and 0's on $C_{s(kt)}\backslash K_i$ and hence $\0_{k-1}$ on rows $C$. 
    
    We will show below that we can choose $D\subset D_{s(kt)}\backslash\cup_{i=1}^{t}D_{s(i)}$ with $|D|=k$. Then we can find $t\cdot F_{0,k,k-1,0}$ as follows. 
    We have one column in $X_{s(i)}$ for each $i=1,2,\ldots ,t$ which is $\1_{k-1}$ on rows $C$ and $\0_{k}$ on rows $D$ (since $D\subset C_{s(i)}\backslash C$ for each 
    $i=1,2,\ldots ,t$).
    The $t$ columns of $X_{s(tk)}$ we have selected have $\0_{k-1}$ on rows $C$ and 1's on $D_{s(kt)}$ where $D\subseteq D_{s(kt)}$ and hence $\1_{k}$ on rows $D$.      
  These $2t$ columns yield $t\cdot F_{0,k,k-1,0}$ in $[X_{s(1)}\,|\,X_{s(2)}\,\cdots |\,X_{s(t)}\,|\,X_{s(kt)}]$.
    
    To show that  $D$ can be chosen  we first show that $D_{s(i)}\backslash D_{s(j)}\le k-2$ for $s(i)<s(j)$.
    Assume the contrary, $D_{s(i)}\backslash D_{s(j)}\ge k-1$ for $s(i)<s(j)$.  We choose $C'\subseteq D_{s(i)}\backslash D_{s(j)}$ with $|C'|=k-1$.
    Given $s(j)> s(i)$, then  $D_{s(j)}\backslash D_{s(i)}\ge k$ and so we may choose $D'\subseteq D_{s(j)}\backslash D_{s(i)}$ with $|D'|=k$. Now $C'\subset C_{s(j)}$ and $D'\subset C_{s(i)}$.  The number of possible columns of $X_{s(j)}$ with at least one 1 on the rows $C'$ is at most $m^{k-2}$ and with $|X_{s(j)}|\ge m^{k-1}+t$, we find $t$ columns of $X_{s(j)}$ with 0's on rows $C'$ and necessarily with 1's on rows $D'$. The number of possible columns of $X_{s(i)}$ with at least one 1 on the rows of $D'$ is $|D'|m^{k-2}<m^{k-1}$. 
    Given $|X_{s(i)}|\ge m^{k-1}+t$, we find $t$ columns of $X_{s(i)}$ with 0's on rows $D'$ and necessarily with 1's on rows $C'$.
    This yields $t\cdot F_{0,k,k-1,0}$ in $[X_{s(i)}\,|\,X_{s(j)}]$, a contradiction.
    Thus $D_{s(i)}\backslash D_{s(j)}\le k-2$ for $s(i)<s(j)$. 
    We may now conclude that $|D_{s(kt}\backslash\cup_{i=1}^{t}D_{s(i)}|\ge k$ and so a choice for $D$ exists. We conclude $\sum_{i\in T(k-1,0)}|X_i|\le (tk)m^{k-1}$. This concludes Case 2.
    
    \vskip 5pt

  There are $k+1$ choices for type $(a,b)$ and so 
  $$\sum_{i=0}^m |X_i|\le \sum_{j=0}^{k}\left(\sum_{i\in T(j,k-1-j)}|X_i|\right)\le (k+1)(2tk)m^{k-1}$$ 
  and so $\ncols{A}\le (2tk(k+1))m^{k-1}+  (6(k-1))^{5k+2}m^{k-2}$ which is $O(m^{k-1})$.
\qed

\vskip 10pt
\noindent{\bf Proof of \trf{chestnut} for $k\ge 3$: } We use \rf{newinductiont}  so that

\noindent $\forb(m,t\cdot F_{0,k,k,0},t-1)\le \forb(m-1,t\cdot F_{0,k,k,0},t-1)+(t-1)\forb(m,\{F_{0,k,k,0}, t\cdot F_{0,k,k-1,0}\}).$

Induction on $m$ and \lrf{smallbd} yields the bound. \qed

\section{Some applications of the Induction}\label{additional}

\begin{lemma}\label{l4.1}
Let $H$ be a given simple matrix satisfying $\forb(m,H)$ is $O(m^{\ell})$. Then
$\forb(m,t\cdot H)$ is $O(m^{\ell+1})$.\label{ttimes}\end{lemma}
\proof We  use the induction \rf{newinductiont}  where $F=t\cdot H$ and $H=\s(F)$. Induction on $m$ yields the desired bound. \qed
\vskip 10pt
Let $K_k$ denote the $k\times 2^k$ of all possible (0,1)-columns on $k$ rows. The following is the fundamental result about forbidden configurations. 

\begin{thm}\label{sauer} [Sauer  \cite{Sauer}, Perles and Shelah
\cite{PS}, Vapnik and Chervonenkis \cite{VC}] We have that
$$\forb(m,K_k)=\binom{m}{k-1}+\binom{m}{k-2}+\cdots +\binom{m}{0}.$$
Thus $\forb(m,K_k)$ is $\Theta(m^{k-1})$.\end{thm}

We can apply this result as follows.

\begin{thm}\cite{F83} Let $F$ be a given $k\times \ell$ (0,1)-matrix. Then $\forb(m,F)$ is $O(m^k)$.
\label{general}\end{thm}
\proof  Let $t$ be the maximum multiplicity of a column in $F$ (of course $t\le\ell$). Then $F\prec t\cdot K_k$ and so $\s(F)\prec K_k$. Now \lrf{ttimes} combined with \trf{sauer} yields the result. \qed

\vskip 10pt
Interestingly this yields the exact result for $\forb(m,2\cdot K_k)$ \cite{G}.   A more precise result of Anstee and F\"uredi \cite{AF} for $\forb(m,t\cdot K_k)$ has the leading term being bounded by $\frac{t+k-1}{k+1}\binom{m}{k}$ for $t\ge  2$.  
 The following surprising result was obtained by Balogh and Bollob\'as.
\begin{thm}\cite{BB} Let $k$ be given. There is a constant $c_k$ with $\forb(m,\{I_k,I_k^c,T_k\})=c_k$.\end{thm}

This yields the following.

\begin{thm}Let $t,k\ge 2$ be given. Then $\forb(m,\{t\cdot I_k,t\cdot I_k^c,t\cdot T_k\})$ is $\Theta(m)$.\end{thm}
\proof Apply \lrf{ttimes}. 
The matrix $I_m\in\Av(m,\{t\cdot I_k,t\cdot I_k^c,t\cdot T_k\})$  shows that $\forb(m,\{t\cdot I_k,t\cdot I_k^c,t\cdot T_k\})$ is $\Theta(m)$. \qed

\vskip 10pt
\lrf{ttimes} is interesting for those $H$ for which $\forb(m,H)$ is $O(m^{\ell})$ and the number of rows in $H$ is bigger than $\ell$ (see \cite{survey} for examples). It is not expected that this will resolve any \emph{boundary cases}, namely those $F$ for which $\forb(m,[F\,|\,\alpha])$ is bigger than $\forb(m,F)$ by a linear factor (or more) for all choices $\alpha$ which are either not present in $F$ or occur at most once in $F$. The previously mentioned $F_6(t)$ and $F_7(t)$ have quite complicated structure and the induction
\rf{newinductiont} does not appear to work directly.

\end{document}